\documentclass[12pt,epsfig]{article}
\usepackage{amsfonts,amssymb,amsmath,comment}
\usepackage{amsthm}
\usepackage{epsfig}
\usepackage{subfig}
\usepackage{float}
\numberwithin{equation}{section}

\newtheorem{Lemma}{Lemma}[section]
\newtheorem{Theorem}[Lemma]{Theorem}
\newtheorem{Proposition}[Lemma]{Proposition}
\newtheorem{Corollary}[Lemma]{Corollary}
\newtheorem{remark}[Lemma]{Remark}
\newtheorem{definition}[Lemma]{Definition}
\newtheorem{example}[Lemma]{Example}

\usepackage{xcolor}
\newtheorem{Fact}[Lemma]{Fact}
\numberwithin{equation}{section}

\def\bt{\begin{Theorem}}
\def\et{\end{Theorem}}
\def\bl{\begin{Lemma}}
\def\el{\end{Lemma}}
\def\bp{\begin{Proposition}}
\def\ep{\end{Proposition}}
\def\bcor{\begin{Corollary}}
\def\ecor{\end{Corollary}}
\def\bpf{\begin{proof}}
\def\epf{\end{proof}}

\def\brem{\begin{remark}\rm}
\def\erem{\hfill$\Box$\end{remark}}
\def\bedef{\begin{definition}}
\def\endef{\end{definition}}

\def\beg{\begin{example}}
\def\eeg{\end{example}}

\def\bef{\begin{Fact}}
\def\eef{\end{Fact}}

\def\bc{\begin{center}}
\def\ec{\end{center}}

\def\vsq{\vskip .25cm}
\def\beq{\begin{equation}}
\def\eeq{\end{equation}}
\def\beqarray{\begin{eqnarray*}}
\def\eeqarray{\end{eqnarray*}}
\def\<{\leftangle}
\def\>{\rightangle}
\def\({\left(}
\def\){\right)}
\def\f{\varphi}
\def\ds{\displaystyle}

\def\<{\langle}
\def\>{\rangle}
\def\q{\quad}

\def\a{\alpha}
\def\b{\beta}
\def\g{\gamma}

\def\d{\delta}
\def\h{\hbox}

\def\t{\tau}

\def\l{\lambda}
\def\e{\varepsilon}

\def\O{\Omega}

\def\w.r.t.{with respect to}
\def\R{{\mathbb{R}}}
\def\N{{\mathbb{N}}}

\def\C{{\mathbb{C}}}

\def\H{{\mathcal H}}

\def\bq{\begin{quote}}
\def\eq{\end{quote}}

\def\bit{\begin{itemize}}
\def\eit{\end{itemize}}
\def\i{\item}
\def\ben{\begin{enumerate}}
\def\een{\end{enumerate}}

\begin{document}

\title{A Regularization for Time-Fractional Backward Heat Conduction Problem with Inhomogeneous Source Function}
\author{Vighnesh V. Alavani,  \;\; P. Danumjaya, \;\; M.Thamban Nair \\
Department of Mathematics, \\BITS Pilani, K.K. Birla Goa Campus, \\Zuarinagar, Goa 403726, INDIA. \\
Email: vighneshalavani@gmail.com; danu@goa.bits-pilani.ac.in; \\ mtnair@goa.bits-pilani.ac.in }
\date{}
\maketitle

\begin{abstract}
\noindent
Recently, Nair and Danumjaya (2023) introduced a new regularization method for the homogeneous time-fractional backward heat conduction problem (TFBHCP) in a one-dimensional space variable, for determining the initial value function. In this paper, the authors extend the analysis done in the above referred paper to a more general setting of an inhomogeneous time-fractional heat equation involving the higher dimensional state variables and a general elliptic operator. We carry out the analysis for the newly introduced regularization method for the  TFBHCP providing optimal order error estimates under a source condition by choosing the regularization parameter appropriately, and also carry out numerical experiments illustrating the theoretical results.\\\\
{\bf Keywords}: Time-fractional backward heat conduction problem, ill-posed problem, regularization, numerical experiments.   \\\\
{\bf AMS Subject Classifications:} 35K57; 35R25; 35R30; 65J20
\end{abstract}
\section{Introduction}
\setcounter{equation}{0}
For the last 30 years, fractional differential equations have been used to model various problems in science and engineering. The most important of these equations is the diffusion process modelled by 
\begin{equation} \label{Eqn1.1}
   \frac{\partial^\alpha u}{\partial t^\alpha} +  Lu = f, \quad \alpha >0, 
\end{equation}
where $L$ is an elliptic operator in an appropriate domain, $u$ denotes the concentration of the diffusing substance and $f$ accounts for any external influences. When  $\alpha=1$, the above equation is reduced to a traditional diffusion equation. Several studies on the diffusion process have suggested that practical diffusion extends beyond the classical or traditional framework \cite{nig}. In subdiffusion ($0 < \alpha < 1$), particles exhibit slower spreading than predicted by classical diffusion models, while in superdiffusion ($1 < \alpha < 2$), the spreading is accelerated. \\

Researchers have found fractional diffusion equations are necessary for capturing various complex transport phenomena. For example, Nigmatulin \cite{nig}  highlighted their effectiveness in modelling universal responses in electromagnetic, acoustic, and mechanical systems. In contrast, Metzler and Klafter \cite{metz} emphasized their use in describing non-Markovian diffusion processes with memory. The fractional diffusion equation model is also appropriate for investigating problems arising in the areas of spatially disordered systems, porous media, fractal media, turbulent fluids and plasmas, biological media with traps, stock price movements, and so on (see \cite{FDE-Pod, Tuan-1}, and the references therein). \\


The forward problem corresponding to (\ref{Eqn1.1}) is as follows: Given the initial condition $u(x,0)=u_0(x)$ and  source function $f$, find the solution $u(\cdot,\cdot)$ satisfying
\begin{eqnarray}
\label{main}
\frac{\partial^{\alpha} u}{\partial t^{\alpha}} + L u(x, t) &=& f(x, t), \quad x \in \Omega,\; t > 0,  \nonumber \\ 
u(x, 0)&=& u_{0}(x), \quad x \in \Omega, \\
u(x, t)& = & 0,  \quad x \in \partial\Omega, \;t > 0,  \nonumber
\end{eqnarray}
where $\Omega$ is a bounded domain in $\R^d$ with smooth boundary $\partial\Omega$ for some $d\in \N$ and the operator $L$ is given by
\begin{equation}\label{eqA}
L v(x) = -\sum_{i=1}^{d} \frac{\partial}{\partial x_{i}}\left(\sum_{j=1}^{d} A_{i j}(x) \frac{\partial}{\partial x_{j}} v(x)\right) + R(x) v(x), \quad x \in \Omega,   
\end{equation}
where the coefficients $R(x) \geq 0$ and $A_{i j}(x)$  are in $L^{\infty}(\Omega)$. The operator $L$ is assumed to be symmetric and uniformly elliptic, that is, the coefficient functions $A_{ij}$ satisfy the conditions $A_{ij}(x)=A_{ji}(x)$ for all $x \in \Omega$ and $i,j \in \{ 1,2, \ldots,d\}$, and there exists $\beta>0$ such that
$$
\sum_{i, j=1}^{d} A_{i j}(x) \xi_{i} \xi_{j} \geq \beta |\xi|^{2}, \quad \forall x \in \bar{\Omega}, \quad \xi \in \R^d.
$$
One of the inverse problems associated with the above forward problem is to determine $u(\cdot,t)$ for $0 \leq t <\tau$ from the knowledge of $u(\cdot,\tau)$ for $\tau>0$. That is to solve the following final value problem:
\begin{eqnarray} \label{EQB1}
\frac{\partial^{\alpha} u}{\partial t^{\alpha}}+ L u(x, t)  &=& f(x, t), \quad x \in \Omega, \; t > 0,\nonumber
\\ u(x, t)&=&0,  \quad x \in \partial\Omega, \; t > 0,  \\ 
u(x, \tau)&=&g(x), \quad x \in \Omega . \nonumber
\end{eqnarray}

It is well-known that the problem of determining $u(\cdot,0)$ from the knowledge of $u(\cdot,\tau)$ is ill-posed. To obtain a stable approximate solution for the ill-posed problem, it is necessary to apply a regularization method. After Sakamoto and Yamamoto's work in \cite{JMAA-Yama} on the time-fractional backward diffusion problem, several papers on regularization were published. In \cite{jun}, the authors Jun-Gang Wang, {\it et al.} have used Tikhonov regularization to solve the time-fractional homogeneous backward diffusion equation, while the authors L. Wang and J. Liu in \cite{wan} employed the mollification method to regularize the problem. In the literature, several other regularization methods have been developed, such as the nonlocal boundary value problem method in \cite{Hao2, Tuan-1} and quasi-reversible methods in \cite{Hao, Liu-Yama}. \\

In \cite{kokila-nair}, Kokila and Nair have observed that the inverse problem of determining $u(\cdot,t)$ for $0<t<\tau$ from $u(\cdot,\tau)$ is well-posed, in the setting of one dimension with $L = -\Delta$.
This idea has been further explored by Nair and Danumjaya in \cite{Nair-danu},  by considering $u(\cdot,t)$ as regularized approximation for $u(\cdot,0)$, with $t>0$ acting as a regularization parameter. In this paper, we extend the consideration in  \cite{Nair-danu}  to an inhomogeneous time fractional backward heat conduction problem involving a more general elliptic operator in higher dimension.\\

An outline of this article is as follows. We discuss some basic definitions, Lemma's, and Propositions for the subsequent use in section 2. Section 3 derives the well-posedness and ill-posedness of the time fractional backward heat conduction problem (TFBHCP). Section 4 deals with the regularization of TFBHCP. The regularization family for the ill-posed inverse problem is introduced, and its convergence is proved. In Section 5, we derive the error estimates for the noisy data in $g$ and source function $f$. Finally, we perform some numerical experiments to validate the theoretical results in Section 6. \\

Throughout this article, $C$ denotes a generic positive constant that may have different values at different places.

\section{Preliminaries}\label{sec-2}
In this section, we introduce some definitions and results that will be useful in the subsequent sections. \\

We recall that the Caputo fractional derivative of a function $\phi:[0,\t] \rightarrow \R$ for $ \alpha \in (0,1)$ is defined by
\begin{eqnarray*}
    \frac{d^{\alpha} \phi}{d t^{\alpha}}=\frac{1}{\Gamma{(1-\alpha})}\int_0^t (t-s)^{-\alpha} \phi'(s) \,d s,  
\end{eqnarray*} 
where $\Gamma(\cdot)$ is the gamma function defined by
$$
\Gamma(s):=\int_0^\infty e^{-t} t^{s-1}dt, \quad s>0.
$$ 
We use the standard notation for the Sobolev spaces and the corresponding norms. In particular, $L^2(\Omega)$ denotes the Hilbert space of square integrable functions on $\Omega$ with its  inner product  and the induced norm denoted by $(\cdot, \cdot)$ and $\| \cdot \|$, respectively.  
Also, for a Banach space $X$, we denote $L^{\infty}(0,\tau;X)$, the space of all functions $\varphi:[0,\t] \rightarrow X$ such that the  map $t \rightarrow \|\varphi(t)\|_X$  belongs to $L^\infty(0, \t)$. 
It is known that  $L^{\infty}(0,\tau;X)$ is a Banach space with norm $
\| \cdot\|_{L^{\infty}}$ defined by 
$$
\| \varphi \|_{L^{\infty}} = \|\f(\cdot)\|_{L^\infty(0, \t; X)} .
$$

%

Next, we consider a proposition which is a consequence of Lax-Milgram theorem (cf. Nair \cite{Nair-fa}). For the completeness of exposition, we provide its detailed proof as well.   

\bp\label{prop-lax}
Let  $\H$  be an infinite dimensional  real Hilbert space with inner product $\<\cdot, \cdot\>_\H$ and let $\H_0$ be a subspace of $\H$ which is a Hilbert space with a stronger inner product $\<\cdot, \cdot\>_{\H_0}$ and such that the embedding of $\H_0$ into $\H$ is a compact operator. 
Let   $B(\cdot, \cdot): \H_0\times \H_0\to \R$  be a  symmetric  and positive bilinear form  on $\H_0$ which is continuous
and coercive, that is,  there exists $c_0>0$ and $\g_0>0$ such that for all $\f, \psi\in \H_0$,
$$B(\f, \psi) \leq c_0\|\f\|_{\H_0}\|\psi\|_{\H_0},$$
$$B(\f, \f) \geq \g_0\|\f\|_{\H_0}^2.$$
 Then,   there exists a sequence $(\l_n)$ of positive real numbers with $\l_n\to \infty $ as $n\to \infty$ and a sequence $(\f_n)$ in $\H_0$ such that 
$$B(\f_n, \psi) = \l_n\<\f_n, \psi\>_\H\q\forall\, \psi\in \H_0,$$
and $\{\f_n: n\in \N\}$ is an orthonormal basis of $\H$.  
\ep

\bpf 
By Lax-Milgram theorem (cf. Nair \cite{Nair-fa}), for every $\xi\in \H$, there exists a unique $\f_\xi\in \H_0$ such that 
$$B(\f_\xi, \psi) = \<\xi, \psi\>_\H \q\forall\, \psi\in \H_0$$
and 
$\|\f_\xi\|_{\H_0}\leq \frac{1}{\g_0}\|\xi\|_\H.$
Let $T_0: \H\to \H_0$ be the map defined by 
$$T_0 \xi = \f_\xi.$$
Then we see that $T_0$ is a bounded linear operator. Since $\H_0$ is compactly embedded in $\H$, the map $T: \H\to \H$ defined by $T\xi = T_0\xi$ for $\xi\in \H$ is a compact operator. Since $B(\cdot, \cdot)$   be a  symmetric  and positive bilinear form,  $T$ is a compact self-adjoint and positive operator. Hence, by the spectral theorem for such operators (cf. Nair \cite {Nair-linop}), the spectrum of $T$ consists of a countably infinite set of positive eigenvalues $\mu_n,\, n\in \N$ with $\mu_n\to 0$ as $n\to \infty$, and corresponding  eigenfunctions $\f_n,\, n\in \N$, form an orthonormal basis for $\H$. Since $\mu_n\f_n = T\f_n \in \H_0$, it follows that $\f_n\in \H_0$ for every $n\in \N$.   Further, taking $\f_n$ in place of $\xi$, we have 
$$\mu_nB(\f_n, \psi) = B(T\f_n, \psi) = \<\f_n, \psi\>_\H\q\forall\, \psi\in \H,\, n\in \N.$$
Thus, we have 
$$B(\f_n, \psi) = \l_n\<\f_n, \psi\>_\H\q\forall\, \psi\in \H_0,$$
where $\l_n:= 1/\mu_n$ for $n\in \N$. Clearly, $\l_n>0$ for every $n\in \N$ and $\l_n\to \infty$ as $n\to \infty$. 
\epf

The sequence $(\l_n)$ and $(\f_n)$ in the above proposition are called the sequences of eigenvalues and the corresponding eigenfunctions of the bilinear form $B(\cdot, \cdot)$.  In the due course, we shall take $\H$ and $\H_0$ to be the spaces $L^2(\O)$ and $H_0^1(\O)$, respectively. In this case, it is known that $H_0^1(\O)$ is compactly embedded in $L^2(\O)$.

\vsq 
Now, we observe that, for a given $h \in L^2(\Omega)$, the equation $Lv=h$ has the weak form 
$$
a(v,w) = (h,w) \;\; \forall w \in H_0^1(\Omega),
$$ 
where $a(\cdot,\cdot)$ is the bilinear form on  $H_0^1(\Omega)$ defined by
$$
a(v, w)=\sum_{i,j=1}^{d}\int_{\Omega}  A_{i j}(x) \frac{\partial v}{\partial x_{j}}  \frac{\partial w}{\partial x_{i}} dx + \int_{\Omega}R(x) v(x) w(x) dx.
$$
Under the  assumptions on $L$ in (\ref{eqA}), we know that $a(\cdot,\cdot)$ is symmetric, continuous and coercive. Hence, by  Proposition \ref{prop-lax},  there exists an increasing sequence $(\l_n)$ of positive real numbers and a sequence   $(\f_n)$  in $H^1_0(\Omega)$ such that   $\{ \f_n,\, n\in \N \}$ is an orthonormal basis for $L^2(\O)$, which  also satisfy
\begin{equation*}
a(\varphi_n,w)=\lambda_n(\varphi_n,w), \;\; \forall w \in H_0^1(\Omega), \; n \in \N. \label{Eqn1.4new}
\end{equation*}

It is known (see \cite{JMAA-Yama}) that for $u_0 \in L^2(\Omega)$ and $f \in L^{\infty}(0,\tau;L^2(\Omega))$, the weak solution for the forward problem (\ref{main}) can be obtained using the eigenvalues $\lambda_n$ and eigenfunctions $\varphi_n$ as 
\begin{equation} 
  u_{\alpha}(\cdot,t)=\sum_{n=1}^{\infty}
\left[E_{\alpha, 1}\left(-\lambda_n t^\alpha\right)\left(u_0,\varphi_n \right) +  F_{\alpha,n}(t) \right]\varphi_n,  \label{Eqn1.4}	
\end{equation}
where  
\begin{equation}
	F_{\alpha,n}(t)=\int_0^t(t-s)^{\alpha-1} E_{\alpha, \alpha}\left(-\lambda_n(t-s)^\alpha\right) \left(f(\cdot,s) ,\varphi_n \right) \,ds  \label{Eqn1.4-F}
\end{equation}
and for $\a>0$ and $\b\in \R$, 
$E_{\a, \b}(\cdot)$ is the {\it Mittag-Leffler function} \cite{HMS} defined by
\begin{equation} 
E_{\a, \b}(z) = \sum_{k=0}^\infty  \frac{z^k}{\Gamma(\a k + \b)}, \q z\in \C. \label{MLF}
\end{equation} 
For the series on the right hand side of  (\ref{Eqn1.4}) to make sense  as an element of $L^2(\O)$, it is necessary that the two series 
\beq\label{sol-series} \sum_{n=1}^{\infty}
|E_{\alpha, 1}\left(-\lambda_n t^\alpha\right)|^2 |\left(u_0,\varphi_n \right)|^2   \q\h{and}\q  \sum_{n=1}^{\infty} |F_{\alpha,n}(t)|^2\eeq
converge. In order to see this, first, we list some properties of $E_{\a,\b}(\cdot)$ in the following known result  (cf. \cite{jin, FDE-KST, FDE-Pod}). 

\bp\label{MLF-1}
For $\a>0$ and $\b\in \R$, let $E_{\a, \b}(\cdot)$ be defined as in (\ref{MLF}). Then, the following results hold.
\ben
\i[\rm(i)]\,Let 
$0 \leq \a \leq 1$ and $\b\geq \a $. Then the function $\xi \rightarrow E_{\a, \b}(\xi)$ is monotonically decreasing. In particular,
$$  
E_{\a,\b}(-\xi) \leq  E_{\a,\b}(0)= \frac{1}{ \Gamma{(\b)}}, \; \forall \; \xi \geq 0.
$$
\i[\rm(ii)]\,  If  $0<\a<2$, $\b\in \R$ and $
\frac{\pi \alpha}{2}< r < \min \{\pi, \pi \a\}
$, then there exists  $C>0$ depending on $(\a, \b, r)$ such that  
\begin{equation*} 
	\left|E_{\alpha, \beta}(z)\right| \leq \frac{C}{1+|z|},  \label{Eqn4.4}
\end{equation*}
for all $z \in \mathbb{C}$ with  $r < |\mbox{arg} (z)| \leq \pi$.
\een 
\ep 
In the due course, we  shall make use of the following particular case of Proposition \ref{MLF-1}(ii).

\bcor\label{MLF-C}
For $0<\a<1$,   there exists $C_\a>0$  such that 
$$	
\left|E_{\alpha, \a+1}(-\xi)\right| \leq \frac{C_\a}{1+\xi }, \; \forall\, \xi>0.
$$
\ecor 

\bpf 
For $0 < \a < 1$ and $\xi>0$, we have $\min  (\pi, \pi \alpha)=\pi \alpha$ and $|\arg (-\xi)|= \pi$. Hence, the conclusion follow  by  Proposition \ref{MLF-1}(ii).
\epf 

By   Proposition \ref{MLF-1}(i),  $ E_{\alpha, 1}\left(-\lambda_n t^\alpha\right)|\leq 1$ for all $n\in \N$ so that the first series in (\ref{sol-series})   converges and 
$$\sum_{n=1}^{\infty}
|E_{\alpha, 1}\left(-\lambda_n t^\alpha\right)|^2 |\left(u_0,\varphi_n \right)|^2\leq \sum_{n=1}^{\infty}
|\left(u_0,\varphi_n \right)|^2 = \|u_0\|^2.$$
To see that the second series in (\ref{sol-series}) converge, we first note that by Proposition \ref{MLF-1}(i),   
\begin{equation*}\label{MLF-2}
	|E_{\alpha, \a}(-\xi)| \leq \frac{1}{\Gamma{(\a)}} \leq 1  \quad \forall \; \xi \geq 0. 
\end{equation*}
Hence,
\beqarray 
|F_{\alpha,n}(t)|  &\leq&  \int_0^t(t-s)^{\alpha-1} |E_{\alpha, \alpha} \left(-\lambda_n(t-s)^\alpha\right)| \, | \left(f(\cdot,s) ,\varphi_n \right)| \,ds \\ 
&\leq&    \int_0^t(t-s)^{\alpha-1}  | \left(f(\cdot,s) ,\varphi_n \right)| \,ds . 
\eeqarray 
Now, using the Cauchy-Schwarz inequality, we have 
\beqarray 
|F_{\alpha,n}(t)|^2 &\leq & \left[  \int_0^t(t-s)^{(\alpha-1)/2}  (t-s)^{(\alpha-1)/2}| \left(f(\cdot,s) ,\varphi_n \right)| \,ds \right]^2 \\
&\leq &   \Big(  \int_0^t(t-s)^{\alpha-1} ds \Big) 
 \Big(  \int_0^t(t-s)^{\alpha-1}  | \left(f(\cdot,s) ,\varphi_n \right)|^2ds \Big) \\ 
&=&  \frac{t^\a}{\a}   
\Big(  \int_0^t(t-s)^{\alpha-1}  | \left(f(\cdot,s) ,\varphi_n \right)|^2ds \Big).
\eeqarray 
Hence, by using Monotone convergence theorem (with respect to counting measure on $\N$) and the fact that   $\|f(\cdot, s)\|\leq \|f\|_{L^\infty}$, we have 
\begin{eqnarray}
\sum_{n=1}^\infty |F_{\alpha,n}(t)|^2 
&\leq &  \frac{t^\a}{\a}  \Big(  \int_0^t(t-s)^{\alpha-1}  \sum_{n=1}^\infty| \left(f(\cdot,s) ,\varphi_n \right)|^2ds \Big) \nonumber \\
& = &  \frac{t^\a}{\a}  \Big(  \int_0^t(t-s)^{\alpha-1}   \|(f(\cdot,s) \|^2 ds \Big) \nonumber \\
&\leq & \|f\|_{L^\infty}^2 \frac{t^{2\a}}{\a^2}. \label{Eqn2.8-F}
\end{eqnarray}
Thus, we have shown that both the series in (\ref{sol-series}) are convergent. \\

Now, we introduce a space that will be used in this paper. For $k\in \N$, let 
\begin{equation*}
    \mathbb{H}^{k}(\Omega)=\Big\{v \in L^{2}(\Omega): \sum_{n=1}^{\infty} \lambda_{n}^{2 k}\left|\left( v, \varphi_{n}\right)\right|^{2}<\infty\Big\}.
\end{equation*}
{It can be easily shown that } $\mathbb{H}^{k}(\Omega)$ is a Hilbert space with the  norm
\begin{equation*}
   \|v\|_{\mathbb{H}^{k}}=\Big(\sum_{n=1}^{\infty} \lambda_{n}^{2 k}\left|\left( v, \varphi_{n}\right)\right|^{2}\Big)^{1 / 2} .
\end{equation*}
Note that  when $k=0$, we have $ \mathbb{H}^0(\Omega)=L^2(\Omega)$. 
\\

Throughout the paper, we shall use the following lemma on estimates for the Mittag-Leffler function defined in (\ref{MLF}). For its proof, we refer to \cite{Liu-Yama, FDE-Pod}.

\bl\label{lem-1} 
Given real numbers $\a_0,\a_1$ such that  $0<\a_0<\a_1<1$, there exists $C_1>0$ and $C_2>0$ such that 
$$ 
\frac{C_1}{\Gamma(1-\a)(1+\chi)} \leq E_{\a,1}(-\chi) \leq   \frac{C_2}{\Gamma(1-\a)(1+\chi)},
$$
for all $\chi>0$ and for all $\a\in [\a_0, \a_1]$. 
\el 

\section{Time-Fractional Backward Heat Conduction Problem}
\setcounter{equation}{0}
In this section, we discuss the well-posedness and ill-posedness of the inverse  problem, {namely, the {\it  time-fractional backward heat conduction problem} (TFBHCP):}

 \begin{quote} 
$P_t$:\, Knowing $g = u_{\a}(\cdot, \tau)$ for some $\tau > 0$ and $f  \in L^{\infty}(0,\tau;L^2(\Omega))$, find $u_{\a}(\cdot, t)$ for $0 \leq t < \tau$.
\end{quote} 

For the  homogeneous TFBHCP in the one-dimensional setting, recently, Nair and Danumjaya in \cite{Nair-danu} have shown that $\{P_t: 0<t<\t\}$ forms a regularization family for obtaining stable approximate solutions for $P_0$. In the present work,  we are extending  analysis in  \cite{Nair-danu}   to the higher dimensional non-homogeneous TFBHCP, resulting in obtaining regularized solutions $u_\a(\cdot,t)$ with  $0<t<\t$ for  $u_\a(\cdot, 0)$. We shall also provide error estimates for $\|u_\a(\cdot,t) - u_\a(\cdot, 0)\|$ under certain {\it a priori} source condition. To our knowledge, no study has been conducted using the above observation, except for the one-dimensional homogeneous setting considered in \cite{Nair-danu}. However, various regularization methods have been discussed recently (see, e.g. \cite{Tuan-Dang, Hao, kokila-nair, Tuan-1}, and the references therein).

\bt\label{exist-unique}
The problem (\ref{EQB1}) has a unique solution if and only if 
\begin{equation} \label{Eqn2.2}
    \sum_{n=1}^{\infty}\left[ \frac{\left(g,\varphi_n \right)- F_{\alpha,n}(\tau)}{(E_{\alpha, 1}\left(-\lambda_n \tau^\alpha\right)}\right]^2 < \infty.
\end{equation}
\et

\bpf 
Suppose the problem (\ref{EQB1}) has a unique solution $u_\a(\cdot,\cdot)$. Then we have  $g = u_{\a}(\cdot, \tau)$ so that   from equation (\ref{Eqn1.4}), we obtain
\begin{equation}
g =\sum_{n=1}^{\infty}
\left(E_{\alpha, 1}\left(-\lambda_n \tau^\alpha\right)\left(u_0,\varphi_n \right) +  F_{\alpha,n}(\tau)\right)\varphi_n. \label{EQNEW3.2}
\end{equation}
Hence,
$$
\left(g,\varphi_n \right) = E_{\alpha, 1}\left(-\lambda_n \tau^\alpha\right)\left(u_0,\varphi_n \right) +  F_{\alpha,n}(\tau), \label{Eqn2.4}
$$
so that 
$$
 \left(u_0,\varphi_n \right) = \frac{\left(g,\varphi_n \right)- F_{\alpha,n}(\tau)}{E_{\alpha, 1}\left(-\lambda_n \tau^\alpha\right)}. \label{Eqn2.5}
$$
Since $u_0\in L^2(\O)$,  that is, $\displaystyle{\sum_{n=1}^\infty} |(u_0, \f_n)|^2<\infty$, we obtain (\ref{Eqn2.2}). 

Conversely, assume that (\ref{Eqn2.2}) holds.  
Let    $v_{0}$  be defined by
\begin{equation*}
   v_0 =\sum_{n=1}^{\infty}\left[ \frac{\left(g,\varphi_n \right)- F_{\alpha,n}(\tau)}{E_{\alpha, 1}\left(-\lambda_n \tau^\alpha\right)}\right]\varphi_n.
\end{equation*}
Then  $v_0\in L^2(\O)$. Hence, from  equation (\ref{main}) with $v_0$ in place of $u_0$ and (\ref{Eqn1.4}), we obtain
\begin{equation*}
u_\alpha(., \tau) = \sum_{n=1}^{\infty} \left(g,\varphi_n \right) \varphi_n = g.
\end{equation*}
This completes the proof. 
\epf  

Throughout this paper, we assume that $g\in L^2(\O)$ satisfies the condition (\ref{Eqn2.2})  in Theorem \ref{exist-unique}, so that the inverse problem under consideration has a unique solution.

\subsection{Issues of well-posedness and  ill-posedness of TFBHCP}
{
In this sub-section, we show that 
the TFBHCP $P_t$  is ill-posed for $t=0$, and it is well-posed for $0<t<\t$.

First we recall from  (\ref{Eqn1.4}) that 
$$
	g=u_{\alpha}(\cdot,\t)=\sum_{n=1}^{\infty}
	\left[E_{\alpha, 1}\left(-\lambda_n \t^\alpha\right)\left(u_0,\varphi_n \right) +  F_{\alpha,n}(\t) \right]\varphi_n.
$$
From this, we have 
$$(g, \f_n) = E_{\alpha, 1}\left(-\lambda_n \t^\alpha\right)\left(u_0,\varphi_n \right) +  F_{\alpha,n}(\t)$$
so that 
$$(u_0,\varphi_n) = 
\frac{\left(g, \varphi_n\right)-F_{\alpha, n}(\tau)}{E_{\alpha, 1}\left(-\lambda_n \tau^\alpha\right)}.$$
Substituting this value of $(u_0,\varphi_n)$ in (\ref{Eqn1.4}), we obtain
\begin{equation}
 u_{\alpha}(\cdot, t)  =\sum_{n=1}^{\infty}\left\{E_{\alpha, 1}\left(-\lambda_n t^\alpha\right)\left[\frac{\left(g, \varphi_n\right)-F_{\alpha, n}(\tau)}{E_{\alpha, 1}\left(-\lambda_n \tau^\alpha\right)}\right]+F_{\alpha, n}(t)\right\} \varphi_n. \label{Eqn3.1}
\end{equation}

Let us assume that the available data is $\tilde g\in L^2(\O)$ in place of  $g$ satisfying Theorem \ref{exist-unique}. 
Replacing $g$ in (\ref{Eqn3.1})  by $\tilde g$, let us formally define  $\tilde{u}_{\alpha}(\cdot, t)$ by 
\begin{equation}
	\tilde{u}_{\alpha}(\cdot, t)  =\sum_{n=1}^{\infty}\left\{E_{\alpha, 1}\left(-\lambda_n t^\alpha\right)\left[\frac{\left(\tilde{g}, \varphi_n\right)-F_{\alpha, n}(\tau)}{E_{\alpha, 1}\left(-\lambda_n \tau^\alpha\right)}\right]+F_{\alpha, n}(t)\right\} \varphi_n. \label{Eqnew2.8}
\end{equation}
It is to be observed that if $\tilde g$ also satisfies the  condition (\ref{Eqn2.2}) in Theorem \ref{exist-unique}, then $\tilde{u}_{\alpha}(\cdot, \cdot)$  is a solution of the TFBHCP with $\tilde g$ in place of $g$. \\

From the equations (\ref{Eqn3.1}) and (\ref{Eqnew2.8}), we obtain
$$u_{\alpha}(\cdot, t) - \tilde{u}_{\alpha}(\cdot, t) = \sum_{n=1}^{\infty}\left[\frac{E_{\alpha, 1}\left(-\lambda_n t^\alpha\right)}{E_{\alpha, 1}\left(-\lambda_n \tau^\alpha\right)}\right]\left(g-\tilde{g}, \varphi_n\right) \varphi_n,$$
so that 
\begin{equation}
	\|u_{\alpha}(\cdot, t) - \tilde{u}_{\alpha}(\cdot, t)\|^2 = \sum_{n = 1}^{\infty} \left| \frac{E_{\alpha, 1}\left(-\lambda_n t^\alpha \right)}{E_{\alpha, 1}\left(-\lambda_n \tau^\alpha\right)} \right|^2  |(g - \tilde{g}, \f_n)|^2. \label{Eqnew2.11}
\end{equation} 
Now, Lemma \ref{lem-1} gives the relations 
\begin{equation}
	\frac{C_1}{C_2} \frac{\left(1+\lambda_n \tau^\alpha\right)}{\left(1+\lambda_n t^\alpha\right)} \leq \frac{E_{\alpha, 1}\left(-\lambda_n t^\alpha\right)}{E_{\alpha, 1}\left(-\lambda_n \tau^\alpha\right)} \leq \frac{C_2}{C_1} \frac{\left(1+\lambda_n \tau^\alpha\right)}{\left(1+\lambda_n t^\alpha\right)}. \label{Eqn3.2}
\end{equation}
Since  $t<\t$, we have 
$t^\a(1+\l_n\t^\a)  < \t^\a(1+\l_nt^\a)$ so that 
\begin{equation*}
  \frac{1+\lambda_n \tau^\alpha}{1+\lambda_n t^\alpha} \leq \frac{\tau^\alpha}{t^\alpha},
\end{equation*}
and from  (\ref{Eqn3.2}) we obtain 
\begin{equation}
\frac{E_{\alpha, 1}\left(-\lambda_n t^\alpha\right)}{E_{\alpha, 1}\left(-\lambda_n \tau^\alpha\right)} \leq   \frac{C_2}{C_1} \left(\frac{\tau}{t}\right)^{\alpha}. \label{Eqn3.3}
\end{equation}
Hence,   (\ref{Eqnew2.11}) implies that 
\begin{equation*}
\|u_{\alpha}(\cdot, t) - \tilde{u}_{\alpha}(\cdot, t)\| \leq  \frac{C_2}{C_1} \left(\frac{\tau}{t}\right)^{\alpha} \|g - \tilde{g}\|. \label{Eqnew2.12}
\end{equation*}
The above inequality shows that the problem of recovering $u_\alpha(\cdot,t)$ for $t>0$ from $g$ is well-posed. In fact, we proved the following theorem.

\bt \label{Th5.1}
Let  $\|g-\tilde{g}\| \leq \delta $ and for $ 0<t<\t$, we have
\begin{equation*}
	\|u_\alpha(\cdot, t)-\tilde{u}_\alpha(\cdot, t)\| \leq \frac{C_2}{C_1}\frac{\tau^\alpha}{t^\alpha}\delta,
\end{equation*}
where $C_1$ and $C_2$ are as in Lemma \ref{lem-1}.
\et

Next, we show that the problem of recovering $u_{\alpha}(\cdot,0)$ from $g$ is ill-posed. 
From (\ref{Eqn1.4-F}) and (\ref{MLF}) we know that  $F_{\a,n}(0) = 0$ and $E_{\a,1}(0)=1$. Hence, at  $t = 0$, equation (\ref{Eqn3.1}) and (\ref{Eqnew2.8}) become
\begin{equation*}
 u_{\alpha}(\cdot, 0) = \sum_{n=1}^{\infty}
 \left[\frac{\left(g, \varphi_n\right)-F_{\alpha, n}(\tau)}{E_{\alpha, 1}\left(-\lambda_n \tau^\alpha\right)}\right] \varphi_n. \label{Eqn3.6} 
 \end{equation*}
and 
\begin{equation*}
\tilde{u}_{\alpha}(\cdot, 0) =  \sum_{n=1}^{\infty}
\left[\frac{\left(\tilde g, \varphi_n\right)-F_{\alpha, n}(\tau)}{E_{\alpha, 1}\left(-\lambda_n \tau^\alpha\right)}\right] \varphi_n, \label{Eqn3.7}
\end{equation*} 
respectively, so that  we arrive at
\begin{equation*}
u_{\alpha}(\cdot, 0)- \Tilde{u}_{\alpha}(\cdot, 0) = \sum_{n=1}^{\infty}\frac{\left(g-\Tilde{g}, \varphi_n\right)\varphi_n}{E_{\alpha, 1}\left(-\lambda_n \tau^\alpha\right)}. 
\end{equation*}
Thus, we have
$$
\|u_\alpha(\cdot,0)-\tilde{u}_\alpha(\cdot,0)\|^2 = \sum_{n=1}^{\infty} \frac{\left[\left(g-\tilde{g},\varphi_n\right)\right]^2}{E_{\alpha,1}\left(-\lambda_n \tau^\alpha\right)^2} .$$
At this point, we  observe from Lemma \ref{lem-1} that 
\begin{equation*}
 \frac{1}{E_{\alpha,1}\left(-\lambda_n \tau^\alpha\right)} \geq 	\Gamma(1-\alpha)\frac{\left(1+\lambda_n \tau^\alpha\right)}{C_2} . \label{Eqn3.4}
\end{equation*}
Hence, we arrive at the  inequality 
\begin{eqnarray}
\|u_\alpha(\cdot,0)-\tilde{u}_\alpha(\cdot,0)\|^2 \geq \sum_{n=1}^{\infty} \frac{\Gamma(1-\alpha)^2\left(1+\lambda_n \tau^\alpha\right)^2}{C_2^2}\left|\left(g-\tilde{g},\varphi_n\right)\right|^2 . \label{Eqn3.9}
\end{eqnarray}
Since $\l_n\to \infty$ as $n\to \infty$, the above inequality shows that small change  in $g$ can lead  to large deviation in the solution $u_\alpha(\cdot,0)$.
To see this explicitly, for  $\d>0$ let  
$$\tilde{g} = g + \delta \varphi_k$$ for some $k\in \N$.  Then $ |\left(g - \tilde{g}, \varphi_k \right)| = \delta$ and $ \left(g - \tilde{g}, \varphi_n \right) = 0$ for all $n \neq k$. Hence, in this case, we have  $\| \tilde{g} - g \| \leq \delta$ whereas   (\ref{Eqn3.9}) gives 
\begin{equation*}
\|u_\alpha(\cdot,0)-\tilde{u}_\alpha(\cdot,0)\| \geq  \delta  \frac{\Gamma(1-\alpha)}{C_2} \left(1+\lambda_k \tau^\alpha\right). \label{Eqnew2.19}
\end{equation*}
This shows that for large $k$, a small noise level $\d$ can lead to a large deviation in the solution. Thus, the problem of recovering $u_\alpha(\cdot,0)$ from the data $g$ is ill-posed.\\
 
The above discussion can be summarized as follows.

\bt
The time-fractional backward heat conduction problem (\ref{EQB1}), namely $P_t$, is well-posed for $0 < t < \tau$ and it is ill-posed at $t = 0$.
\et

\brem 
For $0<t\leq \t$, let 
${\cal R}_{t,\a}: L^2(\Omega) \rightarrow L^2(\Omega)$  be defined by  
\begin{equation*}
	{\cal R}_{t,\a} \psi = \sum_{n=1}^{\infty}\left\{\left[\frac{E_{\alpha,1}\left(-\lambda_n t^\alpha\right)}{E_{\alpha,1}\left(-\lambda_n \tau^\alpha\right)}\right]\left[\left(\psi, \varphi_n\right)- F_{\alpha,n}(\tau)\right] + F_{\alpha, n}(t)\right\}\varphi_n. \label{Eqnew3.1}
\end{equation*}
In view of the expression (\ref{Eqn3.1}), we know that 
$$
{\cal R}_{t,\a}g = u_\a(\cdot, t),\q 0<t\leq \t.
$$ 
We note  by Theorem \ref{Th5.1} that ${\cal R}_{t,\a}: L^2(\Omega) \to L^2(\Omega)$ is a well-defined continuous function.  We shall see that
$$
\{{\cal R}_{t,\a}:  0<t\leq \t\}
$$
is a regularization family for recovering $u_\a(\cdot, 0)$ from noisy data $\tilde g$.  
\erem

}
\section{The Regularization}
\setcounter{equation}{0}

Recall from the last section that the problem of recovering $u_\a(\cdot, 0)$ from the data $g$ is ill-posed, whereas the problem of recovering $ u_\a(\cdot, t)$ from the data $g$ is well-posed.  
Hence, the  following theorem shows that $\{u_\a(\cdot, t): 0<t\leq \t\}$ is a family of regularised solutions for the above mentioned ill-posed problem.

\bt \label{Th4.1}
Let $g = u_{\a}(\cdot, \tau)$. Then 
 \begin{equation*}
\lim_{t\to 0}\|u_\a (\cdot, t)-u_\a(\cdot, 0)\| =0.\label{conv}
 \end{equation*}
\et

\bpf 
By the representation of $u_\a(\cdot, t)$ in (\ref{Eqn1.4}), we have
\begin{eqnarray*}
	u_\alpha(\cdot, t)- u_\a(\cdot, 0) = \sum_{n=1}^{\infty}\left[\left(E_{\alpha,1}\left(-\lambda_n t^\alpha\right)-1\right)\left(u_0, \varphi_n\right)+F_{\alpha,n}(t)\right] \varphi_n.   
\end{eqnarray*}
Then
\begin{equation*}
	\|u_\alpha(\cdot, t)- u_\a(\cdot, 0)\|^2 = \sum_{n=1}^{\infty}\left[\left(E_{\alpha,1}\left(-\lambda_n t^\alpha\right)-1\right)\left(u_0, \varphi_n\right)+F_{\alpha,n}(t)\right]^2. \label{Eqn4.1}  
\end{equation*}
Note that the right hand side of the above equation is 
$\ds \int_\N G_t(n)d\mu(n),$
where $G_t: \N\to \R$ for each $t\in (0, \t)$ is defined by  
\begin{equation*}
	G_t(n)=\left[\left(E_{\alpha,1}\left(-\lambda_n t^\alpha\right)-1\right)\left(u_0, \varphi_n\right)+F_{\alpha,n}(t)\right]^2 
\end{equation*} 
and $\mu$ is the counting measure on $\N$. 
We show  that 
\ben\i[(i)] \,  $G_t(n)\to 0$ as $t\to 0$  for each $n\in \N$  and 
\i[(ii)]\,$|G_t(n)|\leq G(n)$ for each $n\in \N$ for some  $G: \N\to [0, \infty)$  such that  $\displaystyle{\sum_{n=1}^\infty} G(n)$ converges. 
\een 
Once (i) and (ii) are proved, by  dominated convergence theorem (cf. Nair \cite{Nair-MI})   we obtain 
$$
\lim_{t \rightarrow 0}  \|u_\alpha(\cdot, t)- u_\a(\cdot, 0)\|^2 = \lim_{t \rightarrow 0}  \int_\N G_t(n)d\mu(n) = 0,$$
which would complete the proof. \\

Since $E_{\alpha,1}\left(-\lambda_n t^\alpha\right) \rightarrow 1$ and $F_{\alpha,n}(t) \rightarrow 0$ as $t \rightarrow 0$, we have  $G_t(n) \rightarrow 0$ as $t \rightarrow 0$,  for each $ n \in \mathbb{N}$. Thus, (i) above is proved.  To  prove (ii), first we observe that 
\begin{equation*} 
	|G_t(n)|  \leq 2\left[\left(E_{\alpha,1}\left(-\lambda_n t^\alpha\right)-1\right)^2\left(u_0, \varphi_n\right)^2+F_{\alpha,n}(t)^2\right]. 
\end{equation*}
Hence, it is enough to know that  the series 
$$\sum_{n=1}^\infty |E_{\alpha,1}(-\lambda_n t^\alpha) - 1|^2(u_0,\varphi_n)^2 \q\h{and}\q  \sum_{n=1}^\infty |F_{\alpha,n}(t)|^2 $$
are convergent.  We have already seen in Section \ref{sec-2} that 
$\displaystyle{\sum_{n=1}^\infty} |F_{\alpha,n}(t)|^2 $
is convergent.  To see that the other series is also convergent, recall from Proposition \ref{MLF-1}(i) that   $0\leq E_{\alpha, 1}(-x)\leq 1$ for $x>0$ so that 
$$	\left| E_{\alpha,1}\left(-\lambda_n t^\alpha\right)-1\right| \leq 1.$$
Hence, we have  
$$\sum_{n=1}^\infty |E_{\alpha,1}(-\lambda_n t^\alpha) - 1|^2(u_0,\varphi_n)^2 \leq \sum_{n=1}^\infty (u_0,\varphi_n)^2 = \|u_0\|^2. $$
Thus, the proof is complete. \epf

\section{Error Estimates Under Noisy  Data}
\setcounter{equation}{0}

In this section, we discuss the error estimates associated with the regularized approximations when there is noise in the final value function $g$ and also when there is some noise in the source function $f$. 


\subsection{General error estimates and convergence} 

We assume that the available functions  in place of   $f$ and $g$ are   $f^\e \in L^{\infty}(0,\tau;L^2(\Omega))$ and $g^\d \in L^2(\O)$, respectively,  such that 
\beq\label{noisy-data} 
\|f - f^\e\|_{L^{\infty}} \leq \e\q\h{and}\q 
\|g - {g}^\d \| \leq \d,
\eeq  for some $\e>0$ and $\d >0$. Substituting $g^\d$ in place of $g$ and $f^\e$ in place of $f$ in (\ref{Eqn3.1}), we define 
\begin{equation}
	u^{\d,\e}_{\alpha}(\cdot, t)  = \sum_{n=1}^{\infty}\left\{E_{\alpha, 1}\left(-\lambda_n t^\alpha\right)\left[\frac{\left(g^\d, \varphi_n\right)- F^\e_{\alpha, n}(\tau)}{E_{\alpha, 1}\left(-\lambda_n \tau^\alpha\right)}\right] + F^\e_{\alpha, n}(t)\right\} \varphi_n, \label{Eqn5.6}
\end{equation}
where $F^\e_{\alpha, n}(t)$ is defined by 
$$
	F_{\alpha,n}^\e(t)=\int_0^t(t-s)^{\alpha-1} E_{\alpha, \alpha}\left(-\lambda_n(t-s)^\alpha\right) \left(f^\e(\cdot,s) ,\varphi_n \right) \,ds.
$$
Recall that by the assumption on $f^\e$, the expression on the right-hand side of the equation (\ref{Eqn5.6}) is well-defined.

\bt \label{Th5.3}
Assume (\ref{noisy-data}).   Then for $ 0<t<\tau ,$
\begin{eqnarray*}   
	\|u_\alpha(\cdot, t)-u^{\d,\e}_\alpha(\cdot, t)\| \leq \sqrt{3}\eta \frac{C_2}{C_1}\frac{\tau^\alpha}{t^\alpha}\left( \frac{\t^{\a}}{\a}+1 \right)+ \sqrt{3}\eta \frac{\t^{\a}}{\a},  
\end{eqnarray*}
where $\eta = \max \{\d, \e\}$, and the constants $C_1$ and $C_2$ are as in Lemma \ref{lem-1}.
\et
\bpf
From (\ref{Eqn3.1}) and (\ref{Eqn5.6}), we see that 
	$\|u_\alpha(\cdot,t)-u^{\d,\e}_\alpha(\cdot,t)\|^2$ is equal to 
	\begin{eqnarray*}
		\sum_{n=1}^{\infty} \left[\frac{E_{\alpha,1}\left(-\lambda_n t^\alpha\right)}{E_{\alpha,1}\left(-\lambda_n \tau^\alpha\right)} \left[\left(g-g^\d,\varphi_n\right)-\left(F_{\alpha,n}(\tau)-F^\e_{\alpha,n}(\tau)\right)\right] 
		+\left(F_{\alpha,n}(t)-F^\e_{\alpha,n}(t)\right)\right]^2.
	\end{eqnarray*}
Now, using the inequality $(a + b + c)^2 \leq 3 \left(a^2 + b^2 + c^2 \right)$  and (\ref{Eqn3.3}), we arrive at
\begin{eqnarray}
	\|u_\alpha(\cdot,t)-u^{\d,\e}_\alpha(\cdot,t)\|^2 &\leq&  \sum_{n=1}^{\infty}\left(\frac{C_2}{C_1}\frac{T^\alpha}{t^\alpha}\right)^2\left( 3\left(g-g^\d,\varphi_n\right)^2 + 3 \left(F_{\alpha,n}(\t)-F^\e_{\alpha,n}(\t)\right)^2 \right) \nonumber \\ && +  3 \left(F_{\alpha,n}(t)-F^\e_{\alpha,n}(t)\right)^2. \label{Eqn5.3-new} 
\end{eqnarray}
Replacing $f$ by $f-f^\e$ in inequality (\ref{Eqn2.8-F}), we obtain
\beq
 \sum_{n=1}^\infty |F_{\alpha,n}(t)-F^\e_{\alpha,n}(t)|^2 \leq  \|f-f^\e\|_{L^\infty}^2 \frac{t^{2\a}}{\a^2} \leq   \frac{\t^{2\a}}{\a^2}\e^2. \label{Eqn5.4-new}  
\eeq
Substituting (\ref{Eqn5.4-new}) in (\ref{Eqn5.3-new}), we arrive at
\begin{equation*}   
	\|u_\alpha(\cdot, t)-u^{\d,\e}_\alpha(\cdot, t)\|^2 \leq \left(\frac{C_2}{C_1}\frac{\tau^\alpha}{t^\alpha} \right)^2 \left(3 \d^2 + 3  \frac{\t^{2\a}}{\a^2}\e^2 \right) + 3   \frac{\t^{2\a}}{\a^2}\e^2 . \label{Eqn5.7}
\end{equation*}
We use $\eta = \max \{\d, \e \}$ and the fact $\sqrt{a+b+c}\leq \sqrt{a}+\sqrt{b}+ \sqrt{c}$ if $a,b,c\geq 0$, to obtain
\begin{equation*}   
	\|u_\alpha(\cdot, t)-u^{\d,\e}_\alpha(\cdot, t)\| \leq \sqrt{3}\eta\frac{C_2}{C_1}\frac{\tau^\alpha}{t^\alpha} \left(\frac{\t^{\a}}{\a}+1 \right)+ \sqrt{3}\eta  \frac{\t^{\a}}{\a}.  
\end{equation*}
This completes the proof. 
\epf

Next, we will show that this regularization method converges as $\d\to 0$ and $\e\to 0$ by choosing the regularization parameter $t$ appropriately.

\bt \label{Th5.4}
Let $\|g-g^\d\| \leq \d $ and $\|f-f^\e\|_{L^{\infty}} \leq \e $, and choose $\eta = \max \{ \d, \e \}$. Then for $ 0<t<\tau $ and $\displaystyle{t_\eta = \eta^{\frac{1-\gamma}{\alpha}}< \tau }$, we have
\begin{equation*}
	\|u^{\d,\e}_\alpha(\cdot, t_\eta) - u_\alpha(\cdot, 0)\| \rightarrow 0 \quad \text {as  } \eta \rightarrow 0.    
\end{equation*}
\et
\noindent
\bpf
We have 
\begin{equation*}
	\|u^{\d,\e}_\alpha(\cdot, t)-u_\alpha(\cdot, 0)\|  = \|u^{\d,\e}_\alpha(\cdot, t)-u_\alpha(\cdot, t) + u_\alpha(\cdot, t) - u_\alpha(\cdot, 0)\|. \label{Eqn5.9}
\end{equation*}
Then, using triangle inequality, we obtain
\begin{equation*}
	\|u^{\d,\e}_\alpha(\cdot, t)-u_\alpha(\cdot, 0)\|  \leq \|u^{\d,\e}_\alpha(\cdot, t)-u_\alpha(\cdot, t)\| + \|u_\alpha(\cdot, t) - u_\alpha(\cdot, 0)\|. \label{Eqn5.10}
\end{equation*}
By using Theorem \ref{Th5.3}, we arrive that 
\begin{equation} 
	\|u^{\d,\e}_\alpha(\cdot, t)-u_\alpha(\cdot, 0)\| \leq \sqrt{3}\eta\frac{C_2}{C_1}\frac{\tau^\alpha}{t^\alpha} \left( \frac{\t^{\a}}{\a}+1 \right)+ \sqrt{3}\eta  \frac{\t^{\a}}{\a}+ \|u_\alpha(\cdot, t)-u_\alpha(\cdot, 0)\|. \label{Eqn5.11}
\end{equation}
Substitute $t=t_\eta$ in equation (\ref{Eqn5.11}) and choose $t_\eta$ such that
\begin{equation*}
	\frac{\eta}{t_\eta^\alpha} = \eta^\gamma,  \quad \gamma \in (0,1),
\end{equation*}
to obtain
\begin{eqnarray*} 
	\|u^{\d,\e}_\alpha(\cdot, t_\eta)-u_\alpha(\cdot, 0)\| &\leq& \eta^\gamma\left( \sqrt{3} \frac{C_2}{C_1}\tau^\alpha \left( \frac{\t^{\a}}{\a}+1 \right)+ \sqrt{3}\eta^{1-\gamma} \frac{\t^{\a}}{\a}  \right) \nonumber \\ &+& \|u_\alpha(\cdot, t_\eta)-u_\alpha(\cdot, 0)\|. \label{Eqn5.12}
\end{eqnarray*}
As $\eta  \rightarrow 0$ implies $t_\eta  \rightarrow 0$ and using Theorem \ref{Th4.1},  we obtain
\begin{equation*}
	\|u^{\d,\e}_\alpha(\cdot, t_\eta)-u_\alpha(\cdot, 0)\| \rightarrow 0, \quad \mbox{as}  \quad \eta \rightarrow 0.    
\end{equation*}
This completes the proof. 
\epf

\subsection{Error estimates under source conditions}

It is well-known in the literature on ill-posed problems that an error estimate is possible for the regularized solution only by assuming certain additional {\it smoothness conditions} on the unknown solution, which are known as {\it source conditions}.  In the next theorem we obtain such an error estimate under the assumption that 
$$u_0:=u_\a(\cdot, 0) \in  {\mathbb H}^p$$
for some $p\in (0, 1]$, 
when the data is noise free.

 \bt  \label{Th5.5}
Let $0<t<\tau$ and $u_0 \in \mathbb{H}^p $ for some $p \in (0,1]$. Then  there exists $C_\a>0$ such that 
\begin{equation*}
	\|u_\alpha(\cdot, t)-u_\alpha(\cdot,0)\|
	\leq  \sqrt{2}t^{\a p}\Big[ C_\a\|u_0\|_{\mathbb{H}^p}+\frac{\t^{\a(1-p)}}{\a}\|f\|_{L^{\infty}}\Big],
\end{equation*}
where $C_\a>0$ is as in Corollary \ref{MLF-C}. 
\et

\bpf 
By the representation of $u_\a(\cdot, t)$ in (\ref{Eqn1.4}), we have 
\begin{equation*}
	\|u_\alpha(\cdot, t)-u_\alpha(\cdot,0)\|^2 = \sum_{n=1}^{\infty}\left[\left(E_{\alpha,1}\left(-\lambda_n t^\alpha\right)-1\right)\left(u_0, \varphi_n\right)^2+F_{\alpha,n}(t)\right]^2.  
\end{equation*}
Hence,
\begin{equation}
	\|u_\alpha(\cdot, t)-u_\alpha(\cdot,0)\|^2
	\leq\sum_{n=1}^{\infty}2\left[\left(E_{\alpha,1}\left(-\lambda_n t^\alpha\right)-1\right)^2\left(u_0, \varphi_n\right)^2+F_{\alpha,n}(t)^2\right].  \label{Eqn5.9New}
\end{equation}
We first observe from the definition of $E_{\a,1}(\cdot)$ that 
\begin{eqnarray*}
 E_{\alpha, 1}(z)-1  
=z E_{\alpha, \alpha+1}(z), \; \forall \; z\in \C.\label{Eqn5.10New}  
\end{eqnarray*}
Now, by Corollary \ref{MLF-C},  there exits $C_\a>0$  such that 
\begin{equation*}
 \left|E_{\alpha, \alpha+1}\left(-\lambda_n t^\alpha\right)\right| \leq \frac{C_\a}{1+\lambda_n t^\alpha}. 
\end{equation*}
Thus, 
\begin{equation} 
  \left| E_{\alpha,1}\left(-\lambda_n t^\alpha\right)-1\right|=\left|-\lambda_n t^\alpha E_{\alpha, \alpha+1}\left(-\lambda_n t^\alpha\right)\right| \leq \frac{C_\a\lambda_n t^\alpha}{1+\lambda_n t^\alpha}.   \label{Eqn5.11New}  
\end{equation}
Substituting (\ref{Eqn5.11New}) in (\ref{Eqn5.9New}), we obtain
\begin{equation}
	\|u_\alpha(\cdot, t)-u_\alpha(\cdot,0)\|^2
	\leq  \sum_{n=1}^{\infty}2\Big[\left(\frac{C_\a \lambda_n t^\alpha}{1+\lambda_n t^\alpha}\right)^2\left(u_0, \varphi_n\right)^2+F_{\alpha,n}(t)^2\Big]. \label{Eqn-5.11}   
\end{equation}
We multiply and divide by $\l_n^{2p}$ in (\ref{Eqn-5.11}), we arrive at
\begin{eqnarray}
	\|u_\alpha(\cdot, t)-u_\alpha(\cdot,0)\|^2
	\leq  \sum_{n=1}^{\infty}2\Big[ C_\a^2\left(\frac{ \lambda_n^{1-p} t^\alpha}{1+\lambda_n t^\alpha}\right)^2 \l_n^{2p}\left(u_0, \varphi_n\right)^2+F_{\alpha,n}(t)^2\Big]. \label{Eqn-5.12}
 \end{eqnarray}
 We observe that 
$$ \frac{ \lambda_n^{1-p} t^\alpha}{1+\lambda_n t^\alpha}   = t^{\a p}  \frac{ (\lambda_nt^\a)^{1-p}}{1+\lambda_n t^\alpha} \leq t^{\a p}. $$
Therefore, from (\ref{Eqn-5.12}) , we have 
$$	\|u_\alpha(\cdot, t)-u_\alpha(\cdot,0)\|^2
\leq  \sum_{n=1}^{\infty}2\Big[ C_\a^2  t^{2\a p} \l_n^{2p}\left(u_0, \varphi_n\right)^2+F_{\alpha,n}(t)^2\Big].$$
 Thus, by the definition of $\mathbb{H}^{p}$-norm and the inequality (\ref{Eqn2.8-F}), we obtain
 \begin{equation*}
 	\|u_\alpha(\cdot, t)-u_\alpha(\cdot,0)\|^2 \leq 2\left[ C_\a^2 t^{2\a p} \|u_0\|_{\mathbb{H}^p}^2+\frac{t^{2\a}}{\a^2}\|f\|_{L^{\infty}}^2\right].
 \end{equation*}
Thus, we have  
 $$
 	\|u_\alpha(\cdot, t)-u_\alpha(\cdot,0)\| \leq \sqrt{2}t^{\a p}\left[ C_\a \|u_0\|_{\mathbb{H}^p}+\frac{\t^{\a(1-p)}}{\a}\|f\|_{L^{\infty}}\right].
 $$
 This completes the proof. 
 \epf

From the above theorem we deduce the following theorem by an appropriate choice of the regularization parameter $t$.

\bt \label{Th5.6}
Let  $\|g-g^\d\|\leq\d $ and $\|f-f^\e\|_{L^{\infty}}\leq\e $, and let  $\eta=\max \{\d,\e \} $.  
Assume that  $u_0 \in \mathbb{H}^p$ for $p \in (0,1]$ and $\eta$ satisfies $\eta^{\frac{1}{(p+1)\alpha}}<\tau$.  Then,  choosing $t_\eta = \eta^{\frac{1}{(p+1)\alpha}}$,  the inequality   
\begin{equation*} 
	\|u^{\d,\e}_\alpha(\cdot, t_\eta)-u_\alpha(\cdot, 0)\| \leq \widetilde{C} \,\eta^{\frac{p}{p+1}}
\end{equation*}
holds,  where $\widetilde{C}>0$ is such that 
$$
 \widetilde{C}\geq \sqrt{3} \tau^\a\Big[   \frac{C_2}{C_1} \left( \frac{\t^{\a}}{\a}+1 \right)+  \frac{\t^{\a}}{\a}\Big]   
 + \sqrt{2}\left[ C_\a \|u_0\|_{\mathbb{H}^p}+\frac{\t^{\a(1-p)}}{\a}\|f\|_{L^{\infty}}\right].
$$
with  $C_1$ and $C_2$ are as in Lemma \ref{lem-1}, and $C_\a$ is as in Corollary \ref{MLF-C}. 
\et

\bpf 
From 
Theorem \ref{Th5.3}  and  Theorem \ref{Th5.5}, we have 
$$
	\|u^{\d,\e}_\alpha(\cdot, t)-u_\alpha(\cdot, t)\| \leq \sqrt{3}\eta\frac{C_2}{C_1}\frac{\tau^\alpha}{t^\alpha} \left( \frac{\t^{\a}}{\a}+1 \right)+ \sqrt{3}\eta  \frac{\t^{\a}}{\a},$$
and 
 $$
\|u_\alpha(\cdot, t)-u_\alpha(\cdot,0)\| \leq \sqrt{2}t^{\a p}\left[ C_\a \|u_0\|_{\mathbb{H}^p}+\frac{\t^{\a(1-p)}}{\a}\|f\|_{L^{\infty}}\right]
$$
respectively.  Now, choosing $t:=t_\eta$ as  
$$t_\eta=\eta^{\frac{1}{(p+1)\a}}$$
we obtain 
$$
\frac{\eta}{t_\eta^\a}=t_\eta^{\a p}=\eta^{\frac{p}{(p+1)}}.
$$
Therefore, 
\beqarray 
	\|u^{\d,\e}_\alpha(\cdot, t_\eta)-u_\alpha(\cdot, t_\eta)\|  &\leq &   \sqrt{3}\eta^{\frac{p}{p+1}} \tau^\a \frac{C_2}{C_1} \left( \frac{\t^{\a}}{\a}+1 \right)+ \sqrt{3}\eta  \frac{\t^{\a}}{\a}\\
& = & \eta^{\frac{p}{p+1}}\Big[ \sqrt{3} \tau^\a \frac{C_2}{C_1} \left( \frac{\t^{\a}}{\a}+1 \right)+ \sqrt{3}\eta^{\frac{1}{p+1}}  \frac{\t^{\a}}{\a}\Big]\\ 
&\leq & \eta^{\frac{p}{p+1}} \sqrt{3} \tau^\a\Big[   \frac{C_2}{C_1} \left( \frac{\t^{\a}}{\a}+1 \right)+  \frac{\t^{\a}}{\a}\Big]\\ 
\eeqarray 	
and 	
$$\|u_\alpha(\cdot, t_\eta)-u_\alpha(\cdot,0)\|   \leq  \sqrt{2}\eta^{\frac{p}{p+1}}\left[ C_\a \|u_0\|_{\mathbb{H}^p}+\frac{\t^{\a(1-p)}}{\a}\|f\|_{L^{\infty}}\right].$$
Hence, from the relation 
\begin{equation*}
	\|u^{\d,\e}_\alpha(\cdot, t_\eta)-u_\alpha(\cdot, 0)\|  \leq \|u^{\d,\e}_\alpha(\cdot, t_\eta)-u_\alpha(\cdot, t_\eta)\| + \|u_\alpha(\cdot, t_\eta) - u_\alpha(\cdot, 0)\|,
\end{equation*}
we arrive at the required estimate. 
\epf

\section{Numerical Illustrations }
\setcounter{equation}{0}

In this section, we shall consider some numerical examples to illustrate the level of approximation of the regularized solutions when $d=2$ and  $L = -\Delta$.\\\\ 
Taking  $\Omega = (0,\pi)\times(0,\pi)$, we  consider a two-dimensional case of problem (\ref{EQB1}) as follows:
\begin{eqnarray*}
\frac{\partial^{\alpha} u}{\partial t^{\alpha}} - \Delta u &=& f(x,y,t), \quad (x,y) \in \Omega,\;\; t \in (0, \tau] , \\
u(x,y,t)&=&0,  \quad (x,y)\in \partial\Omega, \;\; t \in (0,\tau],  \\ 
u(x,y, \tau) &=& g(x,y), \quad (x,y) \in \Omega.
\end{eqnarray*}
We take  
$$ 
f(x,y,t) = (2-\pi^2)\sin{(x)} \sin{(y)}e^{-\pi^2 t},
$$
and $g(x, y)$ is constructed by taking 
$$
u_0(x,y) = \sin(x)\sin(y),
$$
that is, $g= u_\a(\cdot,\t)$, where $u_\a(\cdot, t)$ is given by (\ref{Eqn1.4}) with $u_0$ as above. 
It is known that the eigenvalues and eigenfunctions of the elliptic operator $-\Delta$  are given by
\begin{equation*}
   \lambda_{m,n}=m^2+n^2 \quad \text{and}
   \quad \varphi_{m,n}=\frac{2}{\pi}\sin (mx) \sin (ny), \; \text{ for } \;   m, n = 1,2,3,\ldots. 
\end{equation*}

\subsection{Numerical Procedure} 

For computational purposes, we take $\tau = 1$ and illustrate our theoretical results established for the time-fractional backward heat conduction problem. Below, we discuss the steps that will be used in the numerical computations.

\begin{enumerate}
\item 
We first discretize the given domain $\Omega = (0,\pi)\times(0,\pi)$ into $N$ equal partitions in $x$ and $y$ directions. 
Let $\rho_x = \{x_i\}_{i=0}^N$ and $\rho_y = \{y_j\}_{j=0}^N$ denote  uniform partitions of $[0,\pi]$ such that $x_i=ih, \; i=0,\ldots, N,$ and  $y_j=jh, \; j=0,\ldots, N,$ where $h=\pi/N$. We set $I_i=[x_{i-1},x_i]$, $I_j=[y_{j-1},y_j]$, where $i, j = 1, 2, \ldots, N$. 
\item Next we compute the approximate solution of $u_\alpha(\cdot,t)$ from (\ref{Eqn3.1}) as follows:
\begin{equation*}
 u_{\alpha}(\cdot, t)  = \sum_{m=1}^{30}\sum_{n=1}^{30}\left\{E_{\alpha, 1}\left(-(m^2+n^2)t^\alpha\right)\left[\frac{\left(g, \varphi_{m,n}\right)-F_{\alpha,m, n}(\tau)}{E_{\alpha, 1}\left(-(m^2+n^2)\right)}\right]+F_{\alpha,m, n}(t)\right\} \varphi_{m,n},
\end{equation*}
where the final value $g$ is obtained by using the equation (\ref{EQNEW3.2}) with 
$$
u_0(x,y) = \sin(x)\sin(y), \quad f(x,y,t) = (2-\pi^2)\sin(x) \sin(y)e^{-\pi^2 t},
$$
and
\begin{eqnarray*}
  g &=& \sum_{m=1}^{30}\sum_{n=1}^{30}
\big(E_{\alpha, 1}\left(-(m^2+n^2) \right)\left(u_0,\varphi_{m,n} \right)  \\
&+& \int_0^1(1-s)^{\alpha-1} E_{\alpha, \alpha}\left(-(m^2+n^2)(1-s)^\alpha\right) \left(f(\cdot,s) ,\varphi_{m,n} \right) \,d s \big)\varphi_{m,n}, 
\end{eqnarray*}
here 
\begin{equation*}
 \left(u_0,\varphi_{m,n}  \right) = \int_0^\pi \int_0^\pi u_0(x,y)\frac{2 \sin{(mx)} \sin{(ny)}}{\pi} \, d x \,d y,  
 \end{equation*}
 \begin{equation*}
 \left(f, \varphi_{m,n}  \right) = \int_0^\pi \int_0^\pi f(x,y,t)\frac{2 \sin{(mx)} \sin{(ny)}}{\pi} \, d x \,d y.   
\end{equation*}
The integrals are evaluated using a $4$-point Gaussian quadrature rule. We take $N = 4$ partitions in both $x$ and $y$ directions to perform the quadrature rule. Similarly, we perform the integration in temporal direction $t$ as well. \\

The function $F_{\a, m, n}(\cdot)$ is computed from (\ref{Eqn1.4-F}) and the Mittag-Leffler function is computed using {\bf ml.m} MATLAB file written by Roberto Garrappa \cite{Rob}. 
\item We introduce some noise in the source function $f$ say $f^\e$ and denote the resulting approximate solution as $u^\e_{\alpha}(x,y, t)$. Then we have
$$
f^\e(x,y,t) = f(x,y,t) + \frac{\e}{2}, 
$$
with some noise level $\e > 0$. Note that $\|f^\e - f\|_{L^\infty} \leq \e $. \\

Now compute the approximate solution of $u^\e_{\alpha}(x,y, t)$ from the below expression:
\begin{equation*}
 u^\e_{\alpha}(\cdot, t)  = \sum_{m=1}^{30}\sum_{n=1}^{30}\left\{E_{\alpha, 1}\left(-(m^2+n^2)t^\alpha\right)\left[\frac{\left(g, \varphi_{m,n}\right)-F^\e_{\alpha,m, n}(\tau)}{E_{\alpha, 1}\left(-(m^2+n^2)\right)}\right] + F^\e_{\alpha,m, n}(t)\right\} \varphi_{m,n}. 
\end{equation*}
\item We now introduce the noise in both source function $f$ and the given data $g$, say $f^\e$ and $g^\d$, respectively and denote the resulting approximate solution as $u^{\d,\e}_{\alpha}(x,y, t)$. Then, we have
$$
f^\e(x,y,t) = f(x,y,t) + \frac{\e}{2}, 
$$
and 
$$
g^\d(x,y) = g(x,y) + \frac{\d}{2}. 
$$
We note that
$$
\|f-f^\e\|_{L^{\infty}} \leq \e, \quad \mbox{and} \quad 
\|g-g^\d\| \leq \d.
$$ 
We compute the approximate solution of $u^{\d,\e}_{\alpha}(x,y, t)$ from the below expression:
\begin{equation*}
 u^{\d,\e}_{\alpha}(\cdot, t )  = \sum_{m=1}^{30}\sum_{n=1}^{30}\left\{E_{\alpha, 1}\left(-(m^2+n^2)t^\alpha\right)\left[\frac{\left(g^\d, \varphi_{m,n}\right)-F^\e_{\alpha,m, n}(\tau)}{E_{\alpha, 1}\left(-(m^2+n^2)\right)}\right] + F^\e_{\alpha,m, n}(t )\right\} \varphi_{m,n}. 
\end{equation*}

\end{enumerate}

\subsection{Results and Discussions}

We calculate the errors $\|u_\alpha(\cdot,t)-u_0\|$ corresponding to $t_i=10^{-(i+1)}$ for $i = 1, 2,\ldots 8$ and different values of $\alpha \in \{0.2, 0.4, 0.6, 0.8\}$. The errors are expressed in Table \ref{tab1}. We observe that when $t$ decreases, the errors also decreases. This validates our theoretical results established in Theorem \ref{Th4.1}. We show the solution profiles of $u_0$ and $u_{\a}(x, y, t_i)$ for $\a = 0.8$ in Figure 1. \\

Next we compute $u^\e_{\alpha}(x,y, t_{\e})$ for several values of $\e$ and $\alpha$. We note that $u_0 \in \mathbb{H}^1(\O)$. Therefore by  Theorem \ref{Th5.6}, we  choose $\displaystyle{t_\e = \e^{\frac{1}{2 \a}}}$, take $\e_i = 10^{-(i+2)}, \; i = 1, 2, \ldots, 7$ and different values of $\alpha \in \{0.2, 0.4, 0.6, 0.8\}$. 
The errors $\|u^\e_{\alpha}(x,y, t_{\e}) - u_0\|$ are expressed in Table \ref{tab2}. We observe that when $\e$ decreases, the errors also decrease. The solution profiles of $u_0$ and $u^{\e}_{\a}(x, y, t_{\e})$ for $\a = 0.8$ are shown in Figure 2. \\

For computational purpose we assume $\delta = \e $ so that $\eta = \max \{\delta, \e \} = \e$. Since $u_0 \in \mathbb{H}^1(\O)$, using  Theorem \ref{Th5.6}, we  choose $\displaystyle{t_\eta = \eta^{\frac{1}{2 \a}}}$ and compute $u^{\d,\e}_{\alpha}(x,y, t_{\eta})$ for  several values of $\eta$ and $\alpha$. Here, we take $\eta_i = 10^{-(i+2)}, i = 1,2, \ldots, 7$
and different values of $\alpha \in \{0.2, 0.4, 0.6, 0.8\}$. The errors $\|u^{\d,\e}_{\alpha}(x,y, t_{\eta}) - u_0\|$ are shown in Table \ref{tab3}. We observe that when $\eta$ decreases, the errors decrease. This validates our theoretical results established in Theorem \ref{Th5.6}. We show the solution profiles of $u_0$ and $u^{\d,\e}_{\alpha}(\cdot, t_{\eta})$ for $\a = 0.8$ in Figure 3.

\begin{figure}[H]
  \centering
  \includegraphics[width=0.9\linewidth,height=0.63\textheight]{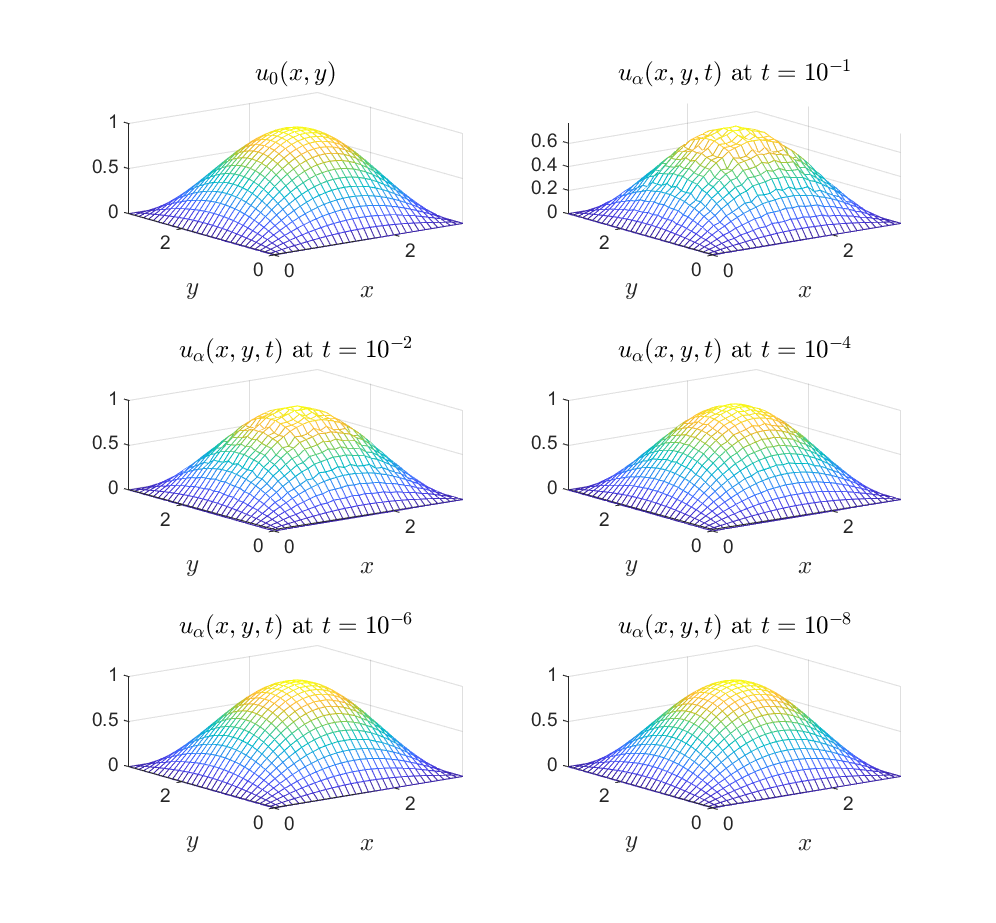}
  \caption{Solution profiles of $u_0(x,y)$ and $u_{\alpha}(x, y, t_{i})$ for $\alpha = 0.8$}
  \label{fig:EX1}
\end{figure}

\begin{table}[H]
\begin{center} 
\caption{$L^2$-errors between $u_0$ and $u_{\alpha}(\cdot, t_i)$}
\label{tab1}
\begin{tabular}{|c|c|c|c|c|} \hline
& $\alpha  = 0.2$ & $\alpha = 0.4$  & $\alpha  = 0.6$ & $\alpha  = 0.8$  \\ \cline{2-5}
$t_i$ \mbox{value} &   $\| u_{\alpha}(\cdot, t)-u_0\|$ & $\| u_{\alpha}(\cdot, t)-u_0 \|$ &  $\|u_{\alpha}(\cdot, t)-u_0\|$  &  $\| u_{\alpha}(\cdot, t) - u_0\|$  \\ \hline  
$10^{-2}$ & $2.3039  $ & $1.8539 $ & $ 9.3104   (-1)$ & $3.8800  (-1)$ \\ \hline
$10^{-3}$ & $1.9435 $ & $9.1668  (-1)$ & $2.6341  (-1)$ & $7.1923  (-2)$ \\ \hline
$10^{-4}$ & $1.4666 $ & $3.9416  (-1)$ & $6.7753  (-2)$ & $1.5789  (-2)$ \\ \hline
$10^{-5}$ & $1.0448  $ & $1.6173  (-1)$ & $1.7109  (-2)$ & $2.8933 (-3)$ \\ \hline
$10^{-6}$ & $7.1524 (-1)$ & $ 6.5167 (-2)$ & $7.0172  (-3)$ & $4.7172 (-4)$ \\ \hline
$10^{-7}$ & $4.7622  (-1)$ & $ 2.6068  (-2)$ & $ 1.8928  (-3)$ & $ 7.0156  (-5)$ \\ \hline
$10^{-8}$ & $3.1115  (-1)$ & $1.0398  (-2)$ &$ 4.8508  (-4)$ & $1.05642  (-5)$ \\ \hline
$10^{-9}$ & 2.0078 $ (-1)$ & 4.1427 $ (-3)$ &  $1.2216 (-4)$ & $ 1.60185  (-6)$ \\ \hline
\end{tabular}
\end{center}
\end{table}

\begin{figure}[H]
\begin{center}
\includegraphics[width=0.9\linewidth, height=0.42\textwidth]{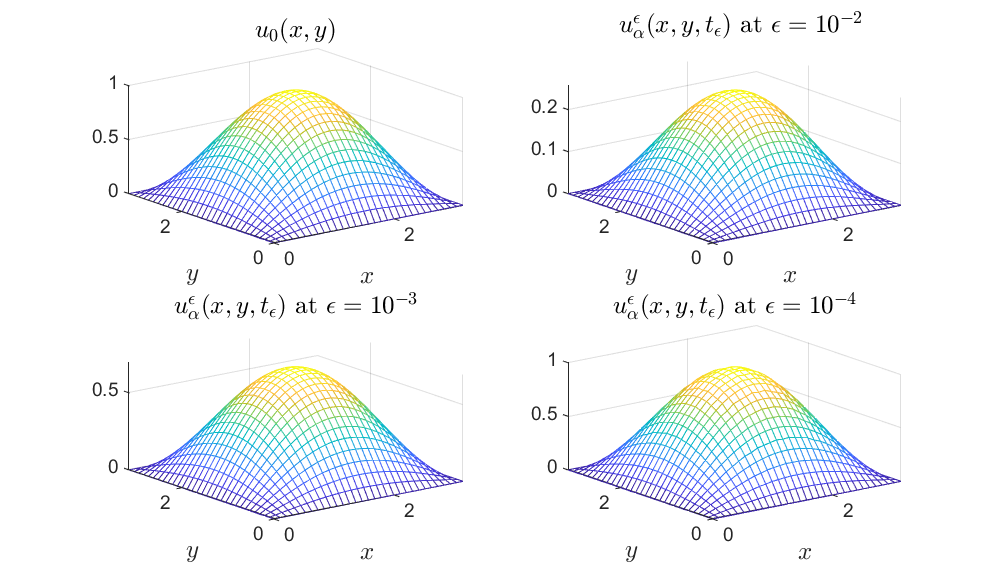}
\caption{Solution profiles of $u_0(x, y)$ and $u^\e_{\alpha}(x,y, t_\e)$ due to noisy in $f$ for $\alpha = 0.8$} 
  \label{fig:EX2}
  \end{center}
\end{figure}

\begin{table}[H] 
\begin{center}
\caption{$L^2$-errors between $u_0$ and $u^\e_{\alpha}(\cdot, t_{\e})$ due to $\e$ noisy in $f$ }
\label{tab2}
\begin{tabular}{|c|c|c|c|c|} \hline
& $\alpha  = 0.2$ & $\alpha = 0.4$  & $\alpha  = 0.6$ & $\alpha = 0.8$  \\ \cline{2-5}
$\e_i$ \mbox{value} &   $\| u^{\e}_{\alpha}(\cdot, t_{\e}) - u_0 \|$ &  $\| u^{\e}_{\alpha}(\cdot, t_{\e}) - u_0 \|$ &  $\| u^{\e}_{\alpha}(\cdot, t_{\e}) - u_0 \|$ &  $\| u^{\e}_{\alpha}(\cdot, t_{\e}) - u_0 \|$  \\ \hline
$10^{-3}$  & $3.8624 (-1)$ & $4.9101 (-1)$ & $5.0883 (-1)$ & $4.7844 (-1)$ \\ \hline
$10^{-4}$ & $1.2858 (-1)$ & $1.6186 (-1)$ & $1.6842 (-1)$ & $1.6200 (-1)$ \\ \hline
$10^{-5}$ & $4.1342 (-2)$ & $5.1862 (-2)$ & $5.3918 (-2)$ & $5.2174 (-2)$ \\ \hline
$10^{-6}$ & $1.5766 (-2)$ & $2.1299 (-2)$ & $2.3100 (-2)$ &$2.3234 (-2)$ \\ \hline
$10^{-7}$  & $6.0274 (-3)$ & $8.0740 (-3)$ & $8.6432 (-3)$  & $8.5687 (-3)$ \\ \hline
$10^{-8}$  & $2.1209 (-3)$ & $2.7968 (-3)$ & $2.9553 (-3)$ & $2.8931 (-3)$\\ \hline
$10^{-9}$  & $6.9930 (-4)$ & $9.1513 (-4)$ & $9.6132 (-4)$ &$9.3612 (-4)$\\ \hline
\end{tabular}
\end{center}
\end{table}

\begin{figure}[H]
\begin{center}
\includegraphics[width=0.9\linewidth, height=0.85\textwidth]{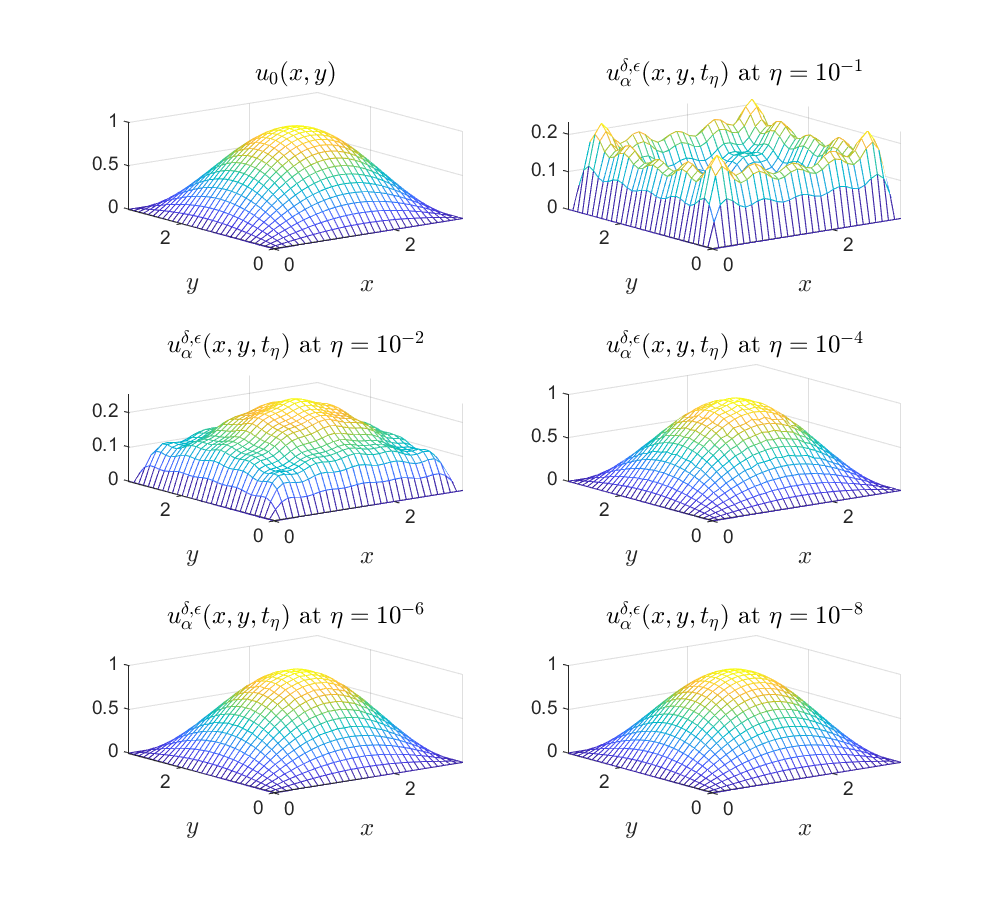}
\caption{Solution profiles of $u_0(x,y)$ and $u^{\d,\e}_{\alpha}(x,y, t_\eta)$  due to noisy in  $f$ and $g$ for $\alpha = 0.8$}.
  \label{fig:EX3}
  \end{center}
\end{figure}
\begin{table}[H]
\begin{center}
\caption{$L^2$-errors between $u_0$ and $u^{\d,\e}_{\alpha}(\cdot, t_{\eta_i})$ due to $\delta$ noisy in $f$ and $g$ }
\label{tab3}
\begin{tabular}{|c|c|c|c|c|} \hline
& $\alpha  = 0.2$ & $\alpha = 0.4$  & $\alpha  = 0.6$ & $\alpha = 0.8$  \\ \cline{2-5}
$\eta_i$ \mbox{value} &   $\| u^{\d,\e}_{\alpha}(\cdot, t_{\eta}) - u_0 \|$ & $\| u^{\d,\e}_{\alpha}(\cdot, t_{\eta}) - u_0 \|$ &  $\| u^{\d,\e}_{\alpha}(\cdot, t_{\eta}) - u_0 \|$ &   $| u^{\d,\e}_{\alpha}(\cdot, t_{\eta}) - u_0 \|$  \\ \hline
  $10^{-3}$ & $3.8253 (-1) $ & $4.8690 (-1) $ & $5.0410 (-1)$ & $4.7320 (-1)$ \\ \hline
$10^{-4}$ & $1.2820 (-1)$ & $1.6142 (-1)$ & $1.6792 (-1)$ & $1.6147 (-1)$ \\ \hline
$10^{-5}$ & $4.1301 (-2)$ & $5.1817 (-2)$ & $5.3867 (-2)$ & $5.2116 (-2)$ \\ \hline
$10^{-6}$ & $1.5764 (-2)$ & $2.1297 (-2)$ & $2.3098 (-2)$ & $2.3234 (-2)$ \\ \hline
$10^{-7}$ & $6.0271 (-3)$ & $8.0737 (-3)$ & $8.6428 (-3)$ & $8.5687 (-3)$ \\ \hline
$10^{-8}$ & $2.1209 (-3)$ & $2.7968 (-3)$ & $2.9553 (-3)$ & $2.8931 (-3)$ \\ \hline
$10^{-9}$ & $6.9930 (-4)$ & $9.1513 (-4)$ & $9.6132 (-4)$ & $9.3612 (-4)$ \\ \hline
\end{tabular}
\end{center}
\end{table}

Below, in Figure 4, we plot the noisy $\eta$ versus error $\|u^{\d,\e}_{\alpha}(\cdot, t_{\eta})-u_0 \|$. We observe that it approximately coincides with $15.2 \sqrt{\eta}$. This validates our theoretical results proved in Theorem \ref{Th5.6}. \\
 
\begin{figure}[H]
 \centering
   \includegraphics[width=0.9\linewidth, height=0.5\textwidth]{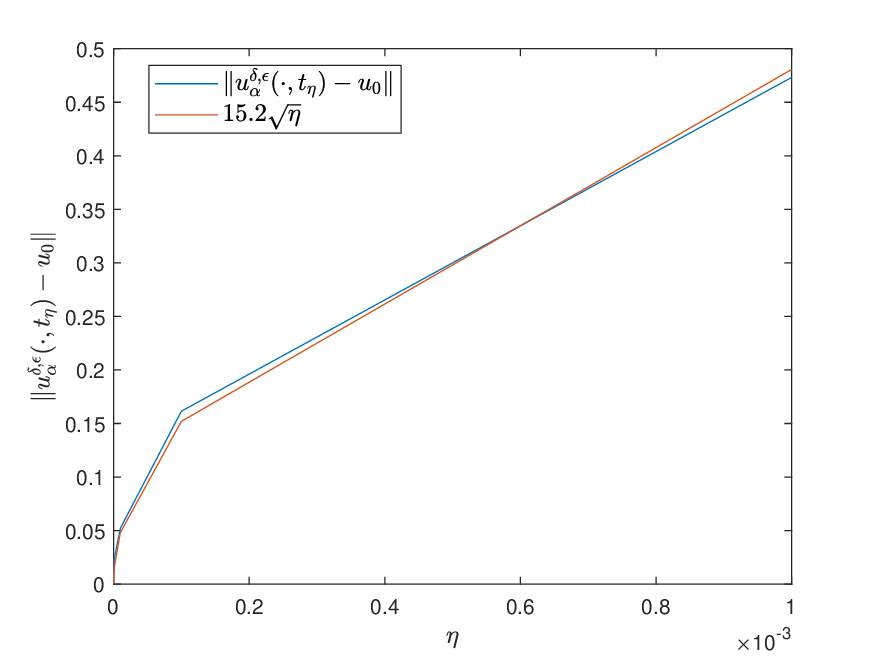}
 \caption{$\eta$ versus error $\|u_{\a}^{\d ,\e}(\cdot,t_{\eta}) - u_0\|$ for $\a = 0.8$. }
  \label{fig:EX4} 
\end{figure}



\end{document}